# ON THE ASYMPTOTIC BEHAVIOR OF COUNTING FUNCTIONS ASSOCIATED TO DEGENERATING HYPERBOLIC RIEMANN SURFACES

Jonathan Huntley
Department of Mathematics
Baruch College CUNY
17 Lexington Avenue
New York, NY
USA

Jay Jorgenson
Department of Mathematics
Yale University
10 Hillhouse Avenue
New Haven, CT 06520
USA

Rolf Lundelius
Department of Mathematics
University of the Witwatersrand
Private Bag 3
P.O. Wits 2050
Republic of South Africa

December 14, 1994

**Abstract.** In this article we will study what we call weighted counting functions on hyperbolic Riemann surfaces of finite volume. If $M$ is compact, then we define the weighted counting function for $w \geq 0$ to be

$$N_{M,w}(T) = \sum_{\lambda_n \leq T} (T - \lambda_n)^w$$

where $\{\lambda_n\}$ is the set of eigenvalues of the Laplacian which acts on the space of smooth functions on $M$. If $M$ is non-compact, then we define the weighted counting function $N_{M,w}(T)$ via the inverse Laplace transform. Using the convergence results from [JL2] concerning the regularized heat trace on finite volume hyperbolic Riemann surfaces, we shall prove the following results.

**Theorem.** *Let $M_\ell$ denote a degenerating family of compact or non-compact hyperbolic Riemann surfaces of finite volume which converges to the non-compact hyperbolic surface $M_0$. For $w \geq 0$ and $T \geq 1/4$, let*

$$G_{\ell,w}(T) = \frac{\Gamma(w+1)}{(16\pi)^{1/2}} \sum_{n=1}^\infty \sum_{\ell_k \in \ell} \frac{\ell_k}{\sinh(n\ell_k/2)} \left(\frac{T - 1/4}{(n\ell_k/2)^2}\right)^{(w+1/2)/2} J_{w+1/2}(n\ell_k\sqrt{(T-1/4)})$$

*where $J_s(x)$ is the J-Bessel function. For $0 \leq T < 1/4$, define $G_{\ell,w}(T)$ to be zero. Then:*
a) *For $w > 3/2$ and $T > 0$, we have*

$$N_{M_\ell,w}(T) = G_{\ell,w}(T) + N_{M_0,w}(T) + o(1) \quad \text{as } \ell \to 0;$$

b) *For $w \geq 0$ and $T > 0$, we have*

$$N_{M_\ell,w}(T) = G_{\ell,w}(T) + o\left(\sum \log(1/\ell_k)\right) \quad \text{as } \ell \to 0;$$

The first author acknowledges support from NSF grant DMS-92-03393, from MSRI under NSF grant DMS-90-22140, and from several PSC-CUNY grants. The second author acknowledges support from NSF grant DMS-93-07023, from the Max-Planck-Institut für Mathematik, and from the Sloan Foundation. Part of the research of the third author was conducted during the US-Sweden Workshop on Spectral Methods sponsored by the NSF under grant INT-9217529.

Typeset by $\mathcal{A}_{\mathcal{M}}\mathcal{S}$-TEX

c) *For $w \geq 0$ and $T \geq 1/4$, we have*
$$G_{\ell,w}(T) = \frac{\Gamma(w+1)(T-1/4)^{w+1/2}}{(4\pi)^{1/2}\Gamma(w+3/2)} \sum_k \log(1/\ell_k) + O(1) \quad \text{as } \ell \to 0.$$

Observe that by taking $T < 1/4$, part (a) proves convergence of certain symmetric functions of the small eigenvalues, thus we prove the main result from [CC]. We discuss a method to improve the error terms in asymptotic expansion of counting functions for $0 \leq w \leq 3/2$.

## §1. Regularized heat traces.

Let $M$ be a connected hyperbolic Riemann surface of finite volume, either compact or non-compact For now, let us assume that $M$ is connected, so then $M$ can be realized as the quotient manifold $\Gamma\backslash\mathbf{h}$, where $\mathbf{h}$ is the hyperbolic plane and $\Gamma$ is a discrete group of isometries of $\mathbf{h}$. Let $K_{\mathbf{h}}(t,\tilde{x},\tilde{y})$ be the heat kernel on $\mathbf{h}$. We shall assume known that the hyperbolic heat kernel on $\mathbf{h}$ is a function of $t > 0$ and of the hyperbolic distance of $\tilde{x}$ to $\tilde{y}$, so
$$K_{\mathbf{h}}(t,\tilde{x},\tilde{y}) = K_{\mathbf{h}}(t,d_{\mathbf{h}}(\tilde{x},\tilde{y})).$$
Quoting from page 246 of [Ch], we have, for $\rho \geq 0$
$$K_{\mathbf{h}}(t,\rho) = \frac{\sqrt{2}e^{-t/4}}{(4\pi t)^{3/2}} \int_\rho^\infty \frac{ue^{-u^2/4t}du}{[\cosh u - \cosh \rho]^{1/2}}$$
and also
$$K_{\mathbf{h}}(t,0) = \frac{1}{2\pi} \int_0^\infty e^{-(1/4+r^2)t} \tanh(\pi r) r\, dr.$$
The heat kernel $K_M(t,x,y)$ on $M$ can be written as a periodization of the heat kernel on the universal cover $\mathbf{h}$, meaning
$$K_M(t,x,y) = \sum_{\gamma \in \Gamma} K_{\mathbf{h}}(t,\tilde{x},\gamma\tilde{y}),$$
where $\tilde{x}$ and $\tilde{y}$ are any points in $\mathbf{h}$ which project via $\Gamma$ to the points $x$ and $y$ in $M$. Let $H(\Gamma)$ denote a set of $\Gamma$-inconjugate primitive hyperbolic classes of $\Gamma$ (meaning classes for which $|\text{Tr}(\gamma)| > 2$ for any representative $\gamma$ of the class), and let $P(\Gamma)$ denote a set of $\Gamma$-inconjugate, non-identity, primitive parabolic classes of $\Gamma$ (meaning classes for which $|\text{Tr}(\gamma)| = 2$ for any representive $\gamma$ of the class). If $M$ is compact, then $P(\Gamma)$ is empty. Let $\Gamma_\gamma$ denote the centralizer of $\gamma$ in $\Gamma$. Recall that the length $\ell(\gamma)$ of the geodesic in the homotopy class determined by $\gamma \in \Gamma$ satisfies the relation $\text{Tr}(\ell(\gamma)) = 2\sinh(\ell(\gamma)/2)$. Now, we can use elementary group theory, as in the derivation of the Selberg trace formula, to write
$$K_M(t,x,y) = K_{\mathbf{h}}(t,\tilde{x},\tilde{y}) + \sum_{n=1}^\infty \sum_{\gamma \in P(\Gamma)} \sum_{\kappa \in \Gamma/\Gamma_\gamma} K_{\mathbf{h}}(t,\tilde{x},\kappa^{-1}\gamma^n\kappa\tilde{y})$$
$$+ \sum_{n=1}^\infty \sum_{\gamma \in H(\Gamma)} \sum_{\kappa \in \Gamma/\Gamma_\gamma} K_{\mathbf{h}}(t,\tilde{x},\kappa^{-1}\gamma^n\kappa\tilde{y})$$
(see [He1], [M], or [S1]). The following theorem defines what we call the hyperbolic heat trace associated to $M$.

**Theorem 1.1.** *Let $M$ be a connected, hyperbolic Riemann surface of finite volume with $p$ cusps. Let $H(\Gamma)$ be a set of inconjugate primitive hyperbolic classes of a uniformizing group $\Gamma$ of $M$.*

a) *For each $t > 0$, the sum*

$$\mathrm{HK}_M(t,x) = \sum_{n=1}^{\infty} \sum_{\gamma \in H(\Gamma)} \sum_{\kappa \in \Gamma/\Gamma_\gamma} K_{\mathbf{h}}(t, \tilde{x}, \kappa^{-1}\gamma^n \kappa \tilde{x})$$

*is a well-defined function on $M$.*

b) *Let $\langle \gamma \rangle$ be the cyclic group generated a representative of $\gamma \in \Gamma$ of a class in $H(\Gamma)$, and let $C_\gamma = \langle \gamma \rangle \backslash \mathbf{h}$. Then*

$$\mathrm{HTr}K_M(t) = \int_M \mathrm{HK}_M(t,x)d\mu(x) = \frac{1}{2} \sum_{\gamma \in H(\Gamma)} \int_{\overset{\circ}{C_\gamma}} \left[ K_{C_\gamma}(t,x,x) - K_{\mathbf{h}}(t,0) \right] d\mu(x).$$

c) *Let $\ell(\gamma)$ be the length of the geodesic in the homotopy class $[\gamma]$ determined by $\gamma \in \Gamma$. Then*

$$\mathrm{HTr}K_M(t) = \frac{e^{-t/4}}{(16\pi t)^{1/2}} \sum_{n=1}^{\infty} \sum_{\gamma \in H(\Gamma)} \frac{\ell(\gamma)}{\sinh(n\ell(\gamma)/2)} e^{-(n\ell(\gamma))^2/4t}.$$

*Outline of Proof.* All aspects of Theorem 1.1 follow from work of Selberg [S2] and calculations in McKean [M] and Hejhal [He1]. Since a detailed discussion is given in [JL1], we simply shall outline the main points. Part (a) follows formal group consideration together with elementary growth considerations for the heat kernel on $\mathbf{h}$ and for the counting function of geodesics on $M$ (see, in particular, Lemma 1.4 of [JL2]). Part (b) follows from unfolding of the integral and integrand, with the factor of $1/2$ appearing since the fundamental group of $C_\gamma$ is isomorphic to $\mathbf{Z}$ whereas the corresponding sum in (a) involves only $\mathbf{Z}^+$ (see also page 233 of [M]). Finally, part (c) follows from the calculations given on page 234 of [M]. □

**Definition 1.2.** *Define the regularized heat trace associated to $M$ by*

$$\mathrm{STr}K_M(t) = \mathrm{HTr}K_M(t) + \mathrm{vol}(M)K_{\mathbf{h}}(t,0).$$

*If $M$ is a hyperbolic Riemann surface of finite volume but not connected, let $M_1, \cdots, M_n$ be the connected components, and define*

$$\mathrm{HTr}K_M(t) = \sum_{j=1}^{n} \mathrm{HTr}K_{M_j}(t) \quad \text{and} \quad \mathrm{STr}K_M(t) = \sum_{j=1}^{n} \mathrm{STr}K_{M_j}(t).$$

Observe that if $M$ is compact, then the regularized trace of the heat kernel is simply the trace of the heat kernel.

In [JL3] we gave a construction of a degenerating family $M_\ell$ of either compact or non-compact hyperbolic Riemann surfaces of finite volume. The reader is referred to this article for complete details, which will be assumed here. Given a $p$-tuple $\ell = (\ell_1, \ell_2, \cdots, \ell_p)$ of positive real numbers, we say that $\ell \to 0$ if the sup norm $|\ell|$ of $\ell$ approaches zero. Although each $M_\ell$ is connected when $\ell > 0$, the limit surface $M_0$ need not be connected, and the number of cusps on $M_0$ is equal to the number of cusps on $M_\ell$ plus $2p$. On a non-connected surface $M$, if $x$ and $y$ lie on different components, $K_M(t,x,y)$ is defined to be zero for all values of $t$.

**Proposition 1.3.** *Let $M_\ell = \Gamma_\ell \backslash \mathbf{h}$ denote a degenerating family of compact or non-compact hyperbolic Riemann surfaces of finite volume which converges to the non-compact hyperbolic surface $M_0$. Let $DH(\Gamma_\ell)$ be the subset of $H(\Gamma_\ell)$ such that the corresponding geodesic has length approaching zero. Define the degenerating heat trace as the integral*

$$\mathrm{DTr} K_{M_\ell}(t) = \frac{1}{2} \sum_{DH(\Gamma_\ell)} \int_{C_{\ell_k}} [K_{C_{\ell_k}}(x,x,t) - K_{\mathbf{h}}(0,t)] d\mu(x).$$

*Then for any $t > 0$ we have the equality*

$$\mathrm{DTr} K_{M_\ell}(t) = \frac{e^{-t/4}}{(16\pi t)^{1/2}} \sum_{n=1}^{\infty} \sum_{\gamma \in DH(\Gamma_\ell)} \frac{\ell(\gamma))}{\sinh(n\ell(\gamma)/2)} e^{-(n\ell(\gamma))^2/4t}.$$

*Proof.* The result follows directly from Theorem 1.1(c). □

**Remark 1.4.** Observe that for every primitive, unoriented geodesic on $M_\ell$ with length going to zero, there are two classes in $H(\Gamma_\ell)$, and the two classes have opposite orientation. Therefore, if we index the sum in Theorem 1.1 by lengths of geodesics rather than elements of $H(\Gamma_\ell)$, the factor of $1/2$ would be eliminated. For another manifestation of this fact, see page 104, both the footnote and line (8.6), of [He1].

The main result of [JL2] is the following convergence theorem.

**Theorem 1.5.** *Let $M_\ell$ denote a degenerating family of compact or non-compact hyperbolic Riemann surfaces of finite volume which converges to the non-compact hyperbolic surface $M_0$. Let $\mathrm{HTr} K_{M_\ell}(t)$ and $\mathrm{DTr} K_{M_\ell}(t)$ be the hyperbolic and degenerating heat traces associated to $M_\ell$.*

a) *(Pointwise) For fixed $t + is$ with $t > 0$, we have*

$$\lim_{\ell \to 0}[\mathrm{STr} K_{M_\ell}(t+is) - \mathrm{DTr} K_{M_\ell}(t+is)] = \mathrm{STr} K_{M_0}(t+is).$$

b) *(Uniformity) For fixed $t > 0$, there is an absolute constant $C$ such that for all $s \in \mathbf{R}$ the bound*

$$|\mathrm{STr} K_{M_\ell}(t+is) - \mathrm{DTr} K_{M_\ell}(t+is)| \leq C(1+|s|^{3/2})$$

*is uniform in $\ell$.*

As an immediate application of Theorem 1.5, we have the following bounds.

**Theorem 1.6.** *For any fixed non-compact, finite volume hyperbolic Riemann surface $M$ and fixed $t > 0$, there is a constant $C = C(M)$ such that*

$$|\mathrm{STr} K_M(t+is)| \leq C(1+|s|)^{3/2}.$$

*Let $M_\ell$ denote a degenerating family of compact or non-compact hyperbolic Riemann surfaces of finite volume which converges to the non-compact hyperbolic surface $M_0$. Then for any fixed $\ell$, there is a constant $C = C(\ell)$ such that*

$$|\mathrm{DTr} K_{M_\ell}(t+is)| \leq C(\ell)(1+|s|)^{3/2}.$$

*Proof.* Line (4.2) in [JL2] states the bound

$$|\text{HTr}K_M(t+is)| \leq C(1+|s|)^{3/2},$$

from which the first assertion follows since $|K_{\mathbf{h}}(z,0)| \leq K_{\mathbf{h}}(t,0)$. The second assertion follows from the first and Theorem 1.5(b). □

## §2. Asymptotics of spectral measures.

For any function $f(t)$ on $\mathbf{R}^+$, we formally define the Laplace transform $\mathcal{L}(f)$ and cumulative distribution function $F$ to be

$$\mathcal{L}(f)(z) = \int_0^\infty e^{-zt} f(t) dt \quad \text{and} \quad F(t) = \int_0^t f(u) du.$$

The Laplace transform $\mathcal{L}(f)$ of $f$ exists if, for example, $f(t)$ is a piecewise continuous, real-valued function for $0 \leq t < \infty$ and there exists and $M$ and $c$ such that $|f(t)| \leq Me^{ct}$. Then $\mathcal{L}(f)(z)$ will exist for some complex $z$ in a half-plane $\text{Re}(z) > a_0$. Recall that the inverse Laplace transform $\mathcal{L}^{-1}$ is

$$f(u) = \frac{1}{2\pi i} \int_{a-i\infty}^{a+i\infty} e^{zt} \mathcal{L}(f)(z) dz \quad \text{and} \quad F(u) = \frac{1}{2\pi i} \int_{a-i\infty}^{a+i\infty} e^{zt} \mathcal{L}(f)(z) \frac{dz}{z}.$$

The inversion formulae holds for any $a > a_0$.

We shall assume that $f$ is such that the Laplace transform $\mathcal{L}(f)$ of $f$ exists, and the inversion formula holds. In addition, we shall make the following basic assumption.

**Assumption:** *There is a constant $a > 0$ such that*

$$\int_{a-i\infty}^{a+i\infty} (1+|\text{Im}(z)|^{3/2}) |\mathcal{L}(f)(z)| \frac{|dz|}{|z|} < \infty.$$

As a direct application of Theorem 1.5, we have the following result.

**Theorem 2.1.** *Let $M_\ell$ denote a degenerating family of compact or non-compact hyperbolic Riemann surfaces of finite volume which converges to the non-compact hyperbolic surface $M_0$. Let $f$ be any function satisfying the stated assumption. Let*

$$N_{M_\ell,S}(f)(T) = \frac{1}{2\pi i} \int_{a-i\infty}^{a+i\infty} \mathcal{L}(f)(z) \text{STr} K_{M_\ell}(z) e^{zT} \frac{dz}{z}$$

*and*

$$N_{M_\ell,D}(f)(T) = \frac{1}{2\pi i} \int_{a-i\infty}^{a+i\infty} \mathcal{L}(f)(z) \text{DTr} K_{M_\ell}(z) e^{zT} \frac{dz}{z}.$$

*Then*

$$\lim_{\ell \to 0} [N_{M_\ell,S}(f)(T) - N_{M_\ell,D}(f)(T)] = N_{M_0,S}(f)(T).$$

*Proof.* The theorem follows from Theorem 1.5 and the dominated convergence theorem. □

## §3. Counting functions for $w > 3/2$.

When considering Theorem 2.1, there is particular interest is the family of functions $f_w(t) = (w+1)t^w$, hence

$$\mathcal{L}(f_w)(z) = \frac{\Gamma(w+2)}{z^{w+1}} \quad \text{and} \quad F_w(t) = t^{w+1}.$$

Let

$$N_{M,w+1}(T) = N_{M,S}((w+1)t^w)(T) = \frac{1}{2\pi i} \int_{a-i\infty}^{a+i\infty} \frac{\Gamma(w+2)}{z^{w+1}} \mathrm{STr} K_M(z) e^{zT} \frac{dz}{z}.$$

We shall call $N_{M,w}(T)$ the *w-th counting function* of the spectrum of the surface $M$. By Theorem 1.6, the weighted counting function is defined for $w > 3/2$. If $M$ is compact, then by Lemma 4.2 below, we can show

$$N_{M,w}(T) = \sum_{\lambda_n \leq T} (T - \lambda_n)^w.$$

Directly from Theorem 2.1 above, we have the following result.

**Theorem 3.1.** *Let $M_\ell$ denote a degenerating family of compact or non-compact hyperbolic Riemann surfaces of finite volume which converges to the non-compact hyperbolic surface $M_0$. For any $w > 3/2$, let*

$$G_{\ell,w}(T) = N_{M_\ell,D}(t^{w-1})(T) = \frac{1}{2\pi i} \int_{a-i\infty}^{a+i\infty} \frac{\Gamma(w+1)\mathrm{DTr} K_M(z)}{z^w} e^{zT} \frac{dz}{z}.$$

*Then for $T > 0$ we have*

$$\lim_{\ell \to 0} [N_{M_\ell,w}(T) - G_{\ell,w}(T)] = N_{M_0,w}(T).$$

Let us now evaluate the asymptotic behavior of the function $G_{\ell,w}(T)$ for fixed $T$ and $w \geq 0$.

**Theorem 3.2.** *Let $M_\ell$ denote a degenerating family of compact or non-compact hyperbolic Riemann surfaces of finite volume.*

a) *For any $w > 3/2$ and $T \geq 1/4$, we have*

$$G_{\ell,w}(T) = \frac{\Gamma(w+1)}{(16\pi)^{1/2}} \sum_{n=1}^{\infty} \sum_{\ell_k \in \ell} \frac{\ell_k}{\sinh(n\ell_k/2)} \left(\frac{T-1/4}{(n\ell_k/2)^2}\right)^{(w+1/2)/2} J_{w+1/2}(n\ell_k\sqrt{(T-1/4)})$$

where $J_s(x)$ is the $J$-Bessel function; if $T < 1/4$, then $G_{\ell,w}(T) = 0$. Therefore, $G_{\ell,w}(T)$ extends to all $w \geq 0$.

b) For fixed $w \geq 0$ and $T \geq 1/4$, we have

$$G_{\ell,w}(T) = \frac{\Gamma(w+1)(T-1/4)^{w+1/2}}{(4\pi)^{1/2}\Gamma(w+3/2)} \sum_k \log(1/\ell_k) + O(1) \quad \text{as } \ell \to 0.$$

*Proof.* To begin, we need two formulas involving the inverse Laplace transform, namely

$$\mathcal{L}^{-1}\left[\mathcal{L}(f)(s)e^{-bs}\right](t) = \begin{cases} f(t-b) & b \leq t < \infty \\ 0 & 0 \leq t < b \end{cases}$$

and, for $\mu > 0$,

$$\mathcal{L}^{-1}\left[s^{-\mu}e^{-k/s}\right](t) = \left(\frac{t}{k}\right)^{(\mu-1)/2} J_{\mu-1}(2\sqrt{(kt)}),$$

where $J$ is the $J$-Bessel function. There are a number of references for these formulas, one of which is [CRC], line 12 on page 465 and line 80 on page 470. Note that the collapse of the $J$-Bessel function, namely

$$J_{-1/2}(x) = \sqrt{\left(\frac{2}{\pi x}\right)} \cos(x)$$

allows for further simplification when $\mu$ is a half-integer, such as stated in the introduction of [JL2]. In any event, these formulas give the equation

$$G_{\ell,w}(T) = \frac{\Gamma(w+1)}{(16\pi)^{1/2}} \sum_{n=1}^{\infty} \sum_{\ell_k \in \ell} \frac{\ell_k}{\sinh(n\ell_k/2)} \left(\frac{T-1/4}{(n\ell_k/2)^2}\right)^{(w+1/2)/2} J_{w+1/2}(n\ell_k\sqrt{(T-1/4)})$$

for $T \geq 1/4$ and $G_{\ell,w}(T) = 0$ for $T < 1/4$, as asserted. FFrom this expression, one immediately obtains the continuation of $G_{\ell,w}(T)$ to any $w \geq 0$.

In order to estimate $G_{\ell,w}(T)$, let us first consider, for each $\ell_k \in \ell$, the sum

$$\frac{\Gamma(w+1)}{(16\pi)^{1/2}} \sum_{n \geq 2/\ell_k} \frac{\ell_k}{\sinh(n\ell_k/2)} \left(\frac{T-1/4}{(n\ell_k/2)^2}\right)^{(w+1/2)/2} J_{w+1/2}(n\ell_k\sqrt{(T-1/4)})$$

Since $n \geq 2/\ell_k$ and $w \geq 0$, we can bound this sum from above by

$$C_T \sum_{n \geq 2/\ell} \frac{\ell_k}{\sinh(n\ell_k/2)} J_{w+3/2}(n\ell_k\sqrt{(T-1/4)}),$$

for some constant $C_T$ which is independent of $\ell_k$. Now let us use the inequality $\sinh(x) \geq e^x/4$ for $x \geq 1$ and the asymptotic formula

$$J_p(x) \sim \sqrt{\left(\frac{2}{\pi x}\right)} \cos(x - \pi/4 - p\pi/2),$$

which we quote from [WW]. Therefore, $J_p(x) \leq Cx^{-1/2} \leq C$ for $x \geq 1$, and the above sum is bounded in absolute value by

$$C_T \sum_{m \geq 2/\ell_k} \ell e^{-m\ell_k/2} \leq C_T \sum_{m=0}^{\infty} \ell e^{-m\ell_k/2} = \frac{C_T \ell_k}{1 - e^{\ell_k/2}},$$

which is bounded independent of $\ell$.

It remains to estimate, for each $\ell_k \in \ell$, the first terms in the series, namely

$$\frac{\Gamma(w+1)}{(16\pi)^{1/2}} \sum_{n \leq 2/\ell_k} \frac{\ell_k}{\sinh(n\ell_k/2)} \left(\frac{T-1/4}{(n\ell_k/2)^2}\right)^{(w+1/2)/2} J_{w+1/2}(n\ell_k \sqrt{(T-1/4)}).$$

For this, we need the estimates

$$\frac{1}{\sinh x} = \frac{1}{x} + O(1/x) \quad \text{and} \quad J_p(x) = \frac{x^p}{\Gamma(p+1)2^p} + O(x^{p+2}) \quad \text{as } x \to 0.$$

Putting these estimates together, we have

$$G_{\ell,w}(T) = \frac{\Gamma(w+1)(T-1/4)^{w+1/2}}{(16\pi)^{1/2}\Gamma(w+1/2)} \sum_k \sum_{0 \leq n \leq 2/\ell_k} \left(\frac{2}{n} + O(n\ell_k^2)\right)\left(1 + O((n\ell_k)^2)\right) + O(1)$$

$$= \frac{\Gamma(w+1)(T-1/4)^{w+1/2}}{(4\pi)^{1/2}\Gamma(w+1/2)} \sum_k \log(1/\ell_k) + O(1),$$

which completes the proof of the theorem. □

One very interesting aspect of Theorem 3.2 is that $G_{\ell,w}(T)$ is zero for $T < 1/4$. Therefore, Theorem 3.2 admits the following corollary.

**Corollary 3.3.** *Let $M_\ell$ denote a degenerating family of compact or non-compact hyperbolic Riemann surfaces of finite volume which converges to the non-compact hyperbolic surface $M_0$. For any $T < 1/4$ and $w > 3/2$, we have*

$$\lim_{\ell \to 0} N_{M_\ell,w}(T) = N_{M_0,w}(T).$$

**Remark 3.4.** From Proposition 4.1 and Proposition 5.1 below, one can write, for any $T < 1/4$,

$$N_{M,w}(T) = \sum_{\lambda_n \leq T} (T - \lambda_n)^w,$$

which means that Corollary 3.3 asserts convergence of certain symmetric functions of the small eigenvalues, which are eigenvalues less than $1/4$. This result has been studied elsewhere, namely [CC] and [He1].

**Remark 3.5.** Let us define the function
$$c_w(T) = \frac{\Gamma(w+1)(T-1/4)^{w+1/2}}{(4\pi)^{1/2}\Gamma(w+3/2)}.$$

Using the identity $s\Gamma(s) = \Gamma(s+1)$, we have
$$\frac{d}{dT}c_{w+1}(T) = (w+1)c_w(T).$$

Observe that for $w \geq 0$, the function $c_w(T)$ is monotone increasing in $T$. These observations will be important in the following section.

## §4. Counting functions for $0 \leq w \leq 3/2$: compact case

In the previous section, we obtained asymptotic behavior of the weighted counting functions for $w > 3/2$. In this section, we consider weighted counted functions for $0 \leq w \leq 3/2$. To begin, we need to improve the bounds given in Theorem 1.6, which is accomplished by using the spectral interpretation of the regularized heat trace as given to us from the Selberg trace formula. FFrom page 44 of [He2], we have the following result.

**Proposition 4.1.** *If $M$ is a compact surface, then the regularized heat trace has the spectral realization*
$$\mathrm{STr}K_M(t) = \sum_{\lambda_n} e^{-\lambda_n t}.$$

From Proposition 4.1, we have that the inverse Laplace transforms which define the weighted counting functions converge for all $w > 0$. Since
$$\frac{d}{dT}N_{M,w+1}(T) = (w+1)N_{M,w}(T)$$
for any $w > 0$, we define $N_{M,0}(T)$ via this differential equation. It is elementary to show
$$N_{M,0}(T) = \sum_{\lambda_n \leq T} 1,$$
where the eigenvalues are counted with multiplicities.

**Lemma 4.2.** *For any function $f$ satisfying the assumptions of Theorem 2.1, we have*
$$\int_0^T F(T-u)\mathcal{L}^{-1}(\mathrm{STr}K)(u)du = \sum_{0 \leq \lambda_n \leq T} F(T-\lambda_n).$$

*Therefore, the distribution $\mathcal{L}^{-1}(\mathrm{STr}K)(u)$ is positive, by which we mean if $f$ is a positive function, the above integral is positive and increasing in $T$.*

*Proof.* Directly from Theorem 8.1 (page 73) and Theorem 12.1(a) (page 91) of [Wi], we have
$$N_{S,M}(f)(T) = \frac{1}{2\pi i}\int_{a-i\infty}^{a+i\infty}\mathcal{L}(f)(z)\mathrm{STr}_M K(z)e^{Tz}\frac{dz}{z} = \int_0^T F(T-u)\mathcal{L}^{-1}(\mathrm{STr}K(t))(u)du$$

and, by Lemma 4.2 and the definition of the Laplace transform,

$$\int_0^T F(T-u)\mathcal{L}^{-1}\left(\text{STr}K(t)\right)(u)du = \sum_{0 \leq \lambda_n \leq T} F(T - \lambda_n),$$

as asserted. Since $f$ is positive, $F$ is positive and increasing, from which the result follows. □

**Lemma 4.3.** *For any $T > 0$ and $w \geq 0$, we have*

$$N_w(T) \leq \frac{1}{w+1}\frac{N_{w+1}(T+\epsilon) - N_{w+1}(T)}{\epsilon} \leq N_w(T+\epsilon).$$

*Proof.* By Lemma 4.2, the function $N_{M,w}(T)$ is a monotone for any $w \geq 0$. With this, both inequalities follow from the mean value theorem. □

The main result of this section is the following theorem.

**Theorem 4.4.** *Let $M_\ell$ denote a degenerating family of compact hyperbolic Riemann surfaces of finite volume which converges to the non-compact hyperbolic surface $M_0$. Then for any $w \geq 0$, and $T \geq 1/4$, we have*

$$N_w(T) \sim c_w(T)\sum \log(1/\ell_k).$$

*Proof.* For fixed $w > 1/2$, $T \geq 1/4$, and $\ell > 0$, consider the function

$$h_{w,\ell}(T) = \frac{N_w(T)}{\sum \log(1/\ell_k)}.$$

From the inequality in Lemma 4.3, we have

$$\frac{1}{w+1}\frac{h_{w+1}(T+\epsilon) - h_{h+1}(T)}{\epsilon} \leq h_{w,\ell}(T+\epsilon).$$

Now let $\ell \to 0$ and use Theorem 3.1 and Theorem 3.2 to get

$$\frac{1}{w+1}\frac{c_{w+1}(T+\epsilon) - c_{w+1}(T)}{\epsilon} \leq \liminf_{\ell \to 0} h_{w,\ell}(T+\epsilon).$$

The right hand side is a decreasing function in $\epsilon$, hence is continuous almost everywhere. Therefore, if we let $\epsilon \to 0$, we get the inequality

$$\frac{1}{w+1}\frac{dc_{w+1}}{dT}(T) \leq \liminf_{\ell \to 0} h_{w,\ell}(T)$$

for almost all $T$. Similarly, we have, for almost all $T$

$$\limsup_{\ell \to 0} h_{w,\ell}(T) \leq \frac{1}{w+1}\frac{dc_{w+1}}{dT}(T).$$

Upon combining, we get the inequality

$$\limsup_{\ell \to 0} h^+_{w,\ell}(T) \leq \frac{1}{w+1} \frac{dc_{w+1}}{dT}(T) \leq \liminf_{\ell \to 0} h_{w,\ell}(T)$$

which holds for almost all $T$. Since the reverse inequalities hold by definition, we have, for almost all $T$,

$$\lim_{\ell \to 0} h_{w,\ell}(T) = \frac{1}{w+1} \frac{dc_{w+1}}{dT}(T) = c_w(T)$$

As observed in Remark 3.5, the function $c_w(T)$ is continuous and monotone, hence the above equality necessarily holds for all $T$.

We now have established the validity of theorem for all $w > 1/2$. We can now repeat the argument for any $w \geq 0$ to finish the proof. □

**Remark 4.5.** For other values of $w \leq 3/2$, the above analysis gives no estimate concerning the error term in the asymptotic expansion of the weighted counting function, unlike the work in [JiZ]. For related remarks, see section 6 below. However, as we shall see in the next section, the above analysis applies to degenerating non-compact hyperbolic Riemann surfaces, a setting not covered by the methods of proof in [He1], [JiZ], or [Wo].

When considering $T < 1/4$ as in Corollary 3.3, the proof of Theorem 4.4 does apply to prove the following result.

**Corollary 4.6.** *Let $M_\ell$ denote a degenerating family of compact hyperbolic Riemann surfaces of finite volume which converges to the non-compact hyperbolic surface $M_0$. If $w > 0$ and $T < 1/4$, then we have the limit*

$$\lim_{\ell \to 0} N_{M_\ell, w}(T) = N_{M_0, w}(T).$$

*If $T < 1/4$ is not an eigenvalue of the Laplacian on the limit surface, then we have the limit*

$$\lim_{\ell \to 0} N_{M_\ell, 0}(T) = N_{M_0, 0}(T).$$

With Corollary 4.6, we obtain many of the results of [CC].

## §5. Counting functions for $0 \leq w \leq 3/2$: non-compact case

The aim of this section is to prove Theorem 4.4 in the case $M_\ell$ is a degenerating family of non-compact surfaces. The method of proof is similar to that of section 4. That is, we shall first quote the Selberg trace formula to get a spectral interpretation of the regularized heat trace. This will allow us to define the spectral counting function for a non-compact surface whenever $w \geq 0$. After this, we shall prove a positivity result associated to the regularized heat trace. After this, the argument of Theorem 4.4 will apply in the non-compact case.

From page 313 of [He3], we have the following spectral interpretation of the regularized heat trace.

**Proposition 5.1.** *If $M$ is a non-compact surface with $p > 0$ cusps, then the regularized heat trace has the spectral realization*

$$\mathrm{STr}K_M(t) = \sum_{C(M)} e^{-\lambda_n t} - \frac{1}{4\pi}\int_{-\infty}^{\infty} e^{-(r^2+1/4)t}\phi'/\phi(1/2+ir)dr$$

$$+ \frac{p}{2\pi}\int_{-\infty}^{\infty} e^{-(r^2+1/4)t}\Gamma'/\Gamma(1+ir)dr - \frac{1}{4}(p-\mathrm{Tr}\Phi(1/2))e^{-t/4} + \frac{p\log 2}{\sqrt{4\pi t}}e^{-t/4},$$

*where $C(M)$ denotes the (possibly finite) set of eigenvalues corresponding to $L^2$ eigenfunctions on $M$, and $\phi(s)$ is the determinant of the scattering matrix $\Phi(s)$*

From Proposition 5.1, we have that the inverse Laplace transforms which define the weighted counting functions converge for all $w > 0$. From the differential equation

$$\frac{d}{dT}N_{M,w+1}(T) = (w+1)N_{M,w}(T),$$

we can define the weighted counting function for all $w \geq 0$. It is an easy exercise, which we leave, to express the counting function $N_{M,0}(T)$ as we define in terms of the spectral data given in Proposition 5.1. We omit the exercise because, as we shall see in the proof of Theorem 5.4, the only factors from Proposition 5.1 which significantly contribute to $N_{M,w}(T)$ are those factors involving the eigenvalues $C(M)$ and the scattering determinant $\phi$.

From pages 72-74 and 160-161 of [He3], specifically line (12.4) and Proposition 12.7, we have the following result.

**Lemma 5.2.** *Let $\{s_k\}$ be the set of numbers for which $1/2 \leq s_k \leq 1$ and $s_k(1-s_k) = \lambda_k \leq 1/4$ is a non-cuspidal eigenvalue of the Laplacian on $M$. Then there is a constant $q_M \geq 0$ such that for all $r$,*

$$-\phi'/\phi(1/2+ir) - \sum_{k=1}^{N}\frac{1-2s_k}{(s_k-1/2)^2+r^2} - 2\log q_M \geq 0.$$

**Remark 5.3.** As on page 160 of [He3], the constant $q_M$ can be described as follows. The scattering determinant can be written as

$$\phi(s) = \left(\frac{\pi^{-(s-1/2)}\Gamma(s-1/2)}{\pi^{-s}\Gamma(s)}\right)^p L(s) \quad \text{where} \quad L(s) = \sum_{\mathbf{q}}\frac{c(\mathbf{q})}{\mathbf{q}^s}$$

for some set of positive numbers $\{\mathbf{q}\}$ tending to infinity. Then by combining the lines (12.3) on page 160, (12.8) on page 166, and the proof and Proposition 12.7 on page 161, we have

$$\log q_M = \inf_{\mathbf{q}}\log\mathbf{q}.$$

As remarked on page 16 of [S2], $\log q_M \leq 0$. Of course, once we have the validity of Lemma 5.2 for one value of $q_M$, the inequality holds when replacing $q_M$ by any smaller, positive number. As we shall see, it is necessary for us to assume that $\log q_M \leq 0$.

We now can argue as in the compact setting of the previous section and obtain the lead term asymptotics of the spectral counting function for degenerating families of non-compact hyperbolic Riemann surfaces.

**Theorem 5.4.** *Let $M_\ell$ denote a degenerating family of non-compact hyperbolic Riemann surfaces of finite volume which converges to the non-compact hyperbolic surface $M_0$. Then for any $w \geq 0$, and $T \geq 1/4$, we have*

$$N_{M_\ell,w}(T) \sim c_w(T) \sum \log(1/\ell_k).$$

*Equivalently, we have the following statement. For any non-compact surface $M$, let*

$$\widehat{\mathrm{STr}K_M}(t) = \sum_{C(M)} e^{-\lambda_n t} - \frac{1}{4\pi} \int_{-\infty}^{\infty} e^{-(r^2+1/4)t} \phi'/\phi(1/2+ir)dr.$$

*Define the distribution*

$$\nu_M(u) = \mathcal{L}^{-1}(\widehat{\mathrm{STr}K_M})(u)$$

*and, for $w \geq 0$, set*

$$\widehat{N_{M,w}}(T) = \int_0^T (T-u)^w \nu_M(u) du.$$

*Then*

$$\widehat{N_{M_\ell,w}}(T) \sim c_w(T) \sum \log(1/\ell_k).$$

*Proof.* In the notation of Proposition 5.1 and Lemma 5.2, let us write the regularized heat trace for any non-compact surface $M$ as

$$\mathrm{STr}K_M(t) = g_{M,1}(t) + g_{M,2}(t) + g_{M,3}(t)$$

where:

$$g_{M,1}(t) = \sum_{C(M)} e^{-\lambda_n t} - \frac{1}{4\pi} \int_{-\infty}^{\infty} e^{-(r^2+1/4)t} \phi'/\phi(1/2+ir)dr$$

$$- \sum_{k=1}^{N} \int_{-\infty}^{\infty} e^{-(r^2+1/4)t} \frac{(1-2s_k)}{4\pi[(s_k-1/2)^2 + r^2]} dr - \frac{\log q_M}{\sqrt{4\pi t}} e^{-t/4};$$

$$g_{M,2}(t) = \sum_{k=1}^{N} \int_{-\infty}^{\infty} e^{-(r^2+1/4)t} \frac{(1-2s_k)dr}{4\pi[(s_k-1/2)^2 + r^2]} + \frac{\log q_M}{\sqrt{4\pi t}} e^{-t/4};$$

$$g_{M,3}(t) = \frac{p}{2\pi} \int_{-\infty}^{\infty} e^{-(r^2+1/4)t} \Gamma'/\Gamma(1+ir)dr - \frac{1}{4}(p - \mathrm{Tr}\Phi(1/2))e^{-t/4} + \frac{p\log 2}{\sqrt{4\pi t}} e^{-t/4}.$$

With this notation, we can write the counting function $N_{M_\ell,w}(T)$ as

$$N_{M_\ell,w}(T) = \int_0^T (T-u)^w \left[ \mathcal{L}^{-1}(g_{M,1})(u) + \mathcal{L}^{-1}(g_{M,2})(u) + \mathcal{L}^{-1}(g_{M,3})(u) \right] du,$$

as well as

$$\widehat{N_{M_\ell,w}}(T) = \int_0^T (T-u)^w \left[ \mathcal{L}^{-1}(g_{M,1})(u) + \mathcal{L}^{-1}(g_{M,2})(u) \right] du,$$

By Lemma 5.2, the function

$$I_{M,1}(T;w) = \int_0^T (T-u)^w \mathcal{L}^{-1}(g_{M,1})(u) du$$

is monotone increasing. Since $\log q_M \leq 0$ and each $s_k \in [1/2, 1]$, the function

$$I_{M,2}(T;w) = \int_0^T (T-u)^w \mathcal{L}^{-1}(g_{M,2})(u) du$$

is monotone decreasing.

Assume $w > 1/2$, $T > 1/4$ is fixed, and let $\epsilon > 0$. By the mean value theorem and the above monotonicity statements, there exist constants $\xi_1, \xi_2 \in [T, T+\epsilon]$ such that

$$I_{M_\ell,1}(T;w) + I_{M_\ell,2}(T+\epsilon;w) + I_{M_\ell,3}(\xi_1;w) \leq \frac{N_{M_\ell,w+1}(T+\epsilon) - N_{M_\ell,w+1}(T)}{\epsilon(w+1)}$$
$$\leq I_{M_\ell,1}(T+\epsilon;w) + I_{M_\ell,2}(T;w) + I_{M_\ell,3}(\xi_2;w).$$

Since $|\text{Tr}(\Phi(1/2))| \leq p$, the function $I_{M_\ell,3}(\xi;w)$ can be bounded independently of $\ell$. Therefore, we have the limit

$$\lim_{\ell \to 0} \frac{I_{M_\ell,3}(\xi;w)}{\sum \log(1/\ell_k)} = 0,$$

hence, we immediately have the limit

$$\lim_{\ell \to 0} \left( \frac{N_{M_\ell,w}(T)}{\sum \log(1/\ell_k)} - \frac{\widehat{N_{M_\ell,w}}(T)}{\sum \log(1/\ell_k)} \right) = 0,$$

which proves the equivalence asserted in the statement of the theorem. With this, we can use Theorem 3.1 and Theorem 3.2 to obtain the inequalities

$$\limsup_{\ell \to 0} \left( \frac{I_{M_\ell,1}(T;w)}{\sum \log(1/\ell_k)} + \frac{I_{M_\ell,2}(T+\epsilon;w)}{\sum \log(1/\ell_k)} \right) \leq \frac{1}{w+1} \frac{c_{w+1}(T+\epsilon) - c_{w+1}(T)}{\epsilon}$$
$$\leq \liminf_{\ell \to 0} \left( \frac{I_{M_\ell,1}(T+\epsilon;w)}{\sum \log(1/\ell_k)} + \frac{I_{M_\ell,2}(T;w)}{\sum \log(1/\ell_k)} \right).$$

If $\epsilon \to 0$ with $\epsilon > 0$, then the lower bound is monotone increasing, and the upper bound is monotone decreasing, since the bounds are limits of monotone functions. Therefore, the upper and lower bounds are continuous almost everywhere; that is, for almost all $T$, we have

$$\limsup_{\ell \to 0} \left( \frac{I_{M_\ell,1}(T;w)}{\sum \log(1/\ell_k)} + \frac{I_{M_\ell,2}(T;w)}{\sum \log(1/\ell_k)} \right) \leq \frac{1}{w+1} \frac{dc_{w+1}}{dT}(T)$$
$$\leq \liminf_{\ell \to 0} \left( \frac{I_{M_\ell,1}(T;w)}{\sum \log(1/\ell_k)} + \frac{I_{M_\ell,2}(T;w)}{\sum \log(1/\ell_k)} \right).$$

Since
$$\widehat{N_{M_\ell,w}}(T) = I_{M_\ell,1}(T;w) + I_{M_\ell,2}(T;w),$$

we have shown
$$\limsup_{\ell\to 0}\left(\frac{\widehat{N_{M_\ell,w}}(T)}{\sum\log(1/\ell_k)}\right) \leq \frac{1}{w+1}\frac{dc_{w+1}}{dT}(T) \leq \liminf_{\ell\to 0}\left(\frac{\widehat{N_{M_\ell,w}}(T)}{\sum\log(1/\ell_k)}\right).$$

The reverse inequalities are trivial. Therefore, we have, for almost all $T$, the existence of the following limit, and the equality
$$\lim_{\ell\to 0}\left(\frac{\widehat{N_{M_\ell,w}}(T)}{\sum\log(1/\ell_k)}\right) = \frac{1}{w+1}\frac{dc_{w+1}}{dT}(T) = c_w(T).$$

Finally, because the function $c_w(T)$ is continuous and monotone increasing in $T$, the equality necessarily holds for all $T$.

Having proved the stated theorem for $w > 1/2$, we can repeat the argument for any $w \geq 0$, which completes the proof of the theorem. $\square$

**Remark 5.5.** By Theorem 5.4, the analogue of Corollary 3.3 and Corollary 4.6 hold in the non-compact setting as well.

## §6. Improved error terms for $0 \leq w \leq 3/2$.

It is seems plausible that one could improve the error term in the asymptotic behavior of the counting functions $N_{M_\ell,w}(T)$ for $w \leq 3/2$ by using arguments as above together with the complication that $\epsilon$ approaches zero at a rate which depends on $\ell$. However, it appears that such results would require one to improve the implied error term in Theorem 1.5(a), which would in turn improve the implied error term in Theorem 3.1. Let us briefly discuss how such results can be obtained.

Assume that for $w > 1/2$ and $T > 1/4$, one has the expansion
$$N_{M_\ell,w+1}(T) = G_{\ell,w+1}(T) + N_{M_0,w+1}(T) + O(f(\ell)) \quad \text{as } \ell \to 0,$$

where $f(\ell)$ is a function which approaches zero as $\ell$ approaches zero. We can argue as in section 4 to get the upper bound
$$(w+1)N_{M_\ell,w}(T) \leq \frac{N_{M_0,w+1}(T+\epsilon) - N_{M_0,w+1}(T)}{\epsilon}$$
$$+ \frac{G_{\ell,w+1}(T+\epsilon) - G_{\ell,w+1}(T)}{\epsilon} + O\left(f(\ell)/\epsilon(\ell)\right).$$

By the mean value theorem, we can write
$$N_{M_\ell,w}(T) \leq N_{M_0,w}(T) + \epsilon(\ell)\frac{d}{dT}N_{M_0,w}(\xi_1) + G_{\ell,w}(T) + \epsilon(\ell)\frac{d}{dT}G_{\ell,w}(\xi_2) + O\left(f(\ell)/\epsilon(\ell)\right)$$

for some $\xi_1, \xi_2 \in [T, T+\epsilon]$. By similar computations, we get
$$N_{M_\ell,w}(T) \geq N_{M_0,w}(T) + \epsilon(\ell)\frac{d}{dT}N_{M_0,w}(\xi_3) + G_{\ell,w}(T) + \epsilon(\ell)\frac{d}{dT}G_{\ell,w}(\xi_4) + O\left(f(\ell)/\epsilon(\ell)\right)$$

for some $\xi_3, \xi_4 \in [T-\epsilon, T]$. Theorem 3.2 now applies. In the end, we have the asymptotic formula

$$N_{M_\ell,w}(T) = N_{M_0,w}(T) + G_{\ell,w}(T) + O\left(\epsilon(\ell)\sum \log(1/\ell_k)\right) + O\left(f(\ell)/\epsilon(\ell)\right).$$

We seek to let $\epsilon = \epsilon(\ell)$ approach zero so that the maximum of the above error terms is minimized. Inductively, one then is able to successively improve the errors for $w > 1/2$ and then for $w \geq 0$. Detailed investigation into the error term in Theorem 1.5(a) will be considered elsewhere.

**Remark 6.1.** In this paper, we considered degenerating hyperbolic Riemann surfaces, so the groups under consideration did not contain elliptic elements. By similar methods, we can address the problem of counting functions associated to the spectral theory of degenerating discrete groups with a fixed set of elliptic elements. In a related question, the spectral asymptotics for degenerating groups with elliptic elements of increasing order, analogous to the results in [JL2], will be presented elsewhere.

*Acknowledgements.* The second author (J. J.) encountered weighted counting functions in collaboration with Serge Lang, based on a course taught by Lang in the academic year 1991-1992 at Yale. At the time, we introduced a technique from analytic number theory, using these counting functions, into our more general context applying also to spectral theory. For instance, we transferred to the general case Theorem 30 on page 82 of [I] and the technique used there. In the present paper, the authors use similar higher weighted counting functions combined with the inverse Laplace transform stemming from work in [W]. The second author very gratefully acknowledges many helpful conversations with Serge Lang concerning the work in the present article and other work. Also, the authors thank Werner Müller whose careful and critical proofreading of [JL2] was one of the factors leading to the considerations in the present article.